\documentclass[a4, 10pt]{amsart}

\usepackage{amssymb}
\usepackage{amstext}
\usepackage{amsmath}
\usepackage{amscd}
\usepackage{latexsym}
\usepackage{amsfonts}
\usepackage{graphicx}
\usepackage[all]{xy}

\theoremstyle{plain}
\newtheorem{thm}{Theorem}[section]
\newtheorem*{thm*}{Theorem}
\newtheorem*{cor*}{Corollary}

\newtheorem{prop}[thm]{Proposition}
\newtheorem{lem}[thm]{Lemma}

\newtheorem{cor}[thm]{Corollary}
\newtheorem{claim}{Claim}
\newtheorem*{claim*}{Claim}

\theoremstyle{definition}
\newtheorem{defn}[thm]{Definition}
\newtheorem{ex}[thm]{Example}
\newtheorem{rem}[thm]{Remark}
\newtheorem{conj}[thm]{Conjecture}
\newtheorem*{conj*}{Conjecture}

\newtheorem{ques}[thm]{Question}

\newtheorem{step}{Step}

\theoremstyle{remark}
\newtheorem*{pf}{{\sl Proof}}

\newtheorem*{ppf}{{\sl Proof of Proposition \ref{pp}}}
\newtheorem*{cpf}{{\sl Proof of Claim}}

\numberwithin{equation}{thm}
\def\Hom{\mathrm{Hom}}

\def\Ext{\mathrm{Ext}}

\def\Gdim{\mathrm{Gdim}}

\def\mod{\mathrm{mod}}

\def\a{\mathfrak a}

\def\m{\mathfrak m}
\def\n{\mathfrak n}

\def\p{\mathfrak p}
\def\q{\mathfrak q}

\def\Z{\Bbb Z}

\def\depth{\mathrm{depth}}

\def\Ass{\mathrm{Ass}}

\def\Min{\mathrm{Min}}

\def\pd{\mathrm{pd}}

\renewcommand{\bar}{\overline}

\begin{document}


\title[On vanishing of certain Ext modules]{On vanishing of certain Ext modules}
\author{Shiro Goto}
\address{Department of Mathematics, School of Science and Technology, Meiji University, 1-1-1, Higashimita, Tama-ku, Kawasaki, Kanagawa 214-8571, Japan}
\email{goto@math.meiji.ac.jp}
\author{Futoshi Hayasaka}
\address{Department of Mathematics, School of Science and Technology, Meiji University, 1-1-1, Higashimita, Tama-ku, Kawasaki, Kanagawa 214-8571, Japan}
\email{fhayasaka@math.meiji.ac.jp}
\author{Ryo Takahashi}
\address{Department of Mathematical Sciences, Faculty of Science, Shinshu University, 3-1-1 Asahi, Matsumoto, Nagano 390-8621, Japan}
\email{takahasi@math.shinshu-u.ac.jp}


\begin{abstract}
Let $R$ be a Noetherian local ring with the maximal ideal $\m$ and $\dim R=1$. 
In this paper, we shall prove that the module $\Ext^1_R (R/Q, R)$ does not vanish for every parameter ideal $Q$ in $R$, 
if the embedding dimension $\mathrm{v}(R)$ of $R$ is at most $4$ 
and the ideal $\m^2$ kills the $0^{\underline{th}}$ local cohomology module $H_{\m}^0(R)$. 
The assertion is no longer true unless $\mathrm{v}(R) \leq 4$. 
Counterexamples are given. We shall also discuss the relation between our counterexamples and a problem on modules of finite G-dimension. 
\end{abstract}

\maketitle
\tableofcontents
\section{Introduction}
Throughout this paper let $R$ be a Noetherian local ring with the maximal ideal $\m$ and $d=\dim R$. 
The purpose of this research is to study the following problem concerning the vanishing of Ext modules. The motivation for the research comes from a conjecture posed by \cite{T} on the modules of 
finite G-dimension.

\begin{ques}\label{ques0}
Let $M$ be an $R$-module of 
finite length. Then does it always hold true that $\Ext_R^d(M,R) \ne (0) $ ?
\end{ques}

In \cite{T} Takahashi studied a characterization of Gorenstein local rings in terms of G-dimension and posed the following conjecture: if a given Noetherian local ring $R$ admits
a non-zero $R$-module of finite length and of finite G-dimension, then the ring $R$ would be Cohen-Macaulay. We can readily see that the conjecture holds true, if Question \ref{ques0} has an affirmative answer. This is the reason why we are interested in Question \ref{ques0}. Later we shall closely discuss the relation between Question \ref{ques0} and the conjecture.

In the present paper we shall restrict our attention on the following  very special case of Question \ref{ques0}. 

\begin{ques}\label{ques}
Assume that $d = 1$ and let $Q$ be a parameter ideal in $R$. Then does it always hold true that 
$\Ext_R^1(R/Q, R) \ne (0)$ ?
\end{ques}

To the surprise of the authors, even in this case the answer is negative in general, while the answer is affirmative in certain special cases even if  the base ring $R$ is not Cohen-Macaulay, as  we shall show in Section 5. Here let us summarize our conclusion into the following two theorems.

\begin{thm}\label{1.3}
Let $d >0$ be an integer. Then there exists a Noetherian local ring 
$R$ such that $\dim R=d$ and $\Ext_R^d(R/Q, R) = (0)$ for some parameter ideal $Q$ in $R$. 
\end{thm}

\begin{thm}\label{main thm}
Let $R$ be a Noetherian local ring with the maximal ideal $\m$ and $\dim R =1$. Assume that  $\frak{m}^2H_{\m}^0(R) = (0)$, where $H_{\m}^0(R)$ denotes the $0^{\underline{th}}$ local cohomology module  of $R$.
Then 
$$ \Ext_R^1(R/Q, R) \ne (0)$$
for every parameter ideal $Q$ in $R$, if $\mathrm{v} (R) \leq 4$. Here $\mathrm{v} (R) = \ell_R(\m / \m^2)$ stands for the embedding dimension of $R$.
\end{thm}

Theorem \ref{main thm} is no longer true unless $\mathrm{v}(R) \leq 4$. 
In Section \ref{count ex} we shall construct examples, which show that for a given integer $v \geq 5$, there exists a parameter ideal $Q$ in a certain one-dimensional 
Noetherian local ring $(R, \m)$ 
with the embedding dimension $\mathrm{v}(R)=v$ and $\m^2 H_{\m}^0(R)=(0)$,  such that 
$\Ext^1_R(R/Q, R)=(0)$. Hence Question \ref{ques} does not hold true in general, and by adding indeterminates to the rings which are one-dimensional counterexamples, we have the negative answer Theorem \ref{1.3} to Question \ref{ques0} for arbitrary dimension $d>0$.

Let us now briefly explain how this paper is organized. We shall prove Theorem \ref{main thm} in Section 3. For the purpose we need some preliminary results and some notation as well, which we will summarize in Section 2.
In Section 4 
we will explore 
some examples affirmative to Question \ref{ques}, which do {\it not} satisfy conditions stated in Theorem \ref{main thm}. In Section 5
we shall prove Theorem \ref{1.3}, constructing counterexamples to Question \ref{ques}. 
In the final Section 6 we will discuss the relation between our counterexamples constructed in Section 5 and the problem on the modules of finite G-dimension. We shall guarantee that the conjecture posed by the third author \cite{T} remains open, showing that our counterexamples given in Section 5 are not  counterexamples for the conjecture of the third author.

\section{Preliminaries}

In this section, we shall summarize some preliminary results which we need to prove Theorem \ref{main thm}. 

Let us fix our notation. Unless otherwise specified, 
let $R$ be a Noetherian local ring with the maximal ideal $\m$ and 
$\dim R = 1$. 
We set $W=H_\m^0(R)$ the $0^{\underline{th}}$ local cohomology module and 
$\a =(0): W$ the annihilator of the ideal $W$. Note that 
$W$ is the unmixed component of $R$, that is, $W=\bigcap_{\p \in \Min R} \q(\p)$, 
where $(0)=\bigcap_{\p \in \Ass R} \q(\p)$ the primary decomposition of $(0)$ in $R$. 
Also, unless otherwise specified, we denote by $Q=(a)$ the parameter ideal in $R$, 
and set $I=(0): Q$. The parameter ideal $Q$ is said to be {\it standard} if $QW=(0)$, that is, $Q$ is contained in $\a$. We denote by $\mu_R(M)$ the minimal number of generators of a 
finitely generated $R$-module $M$, i.e., $\mu_R(M)=\dim_{R/\m}(M/\m M)$. 
We denote by $\mathrm{v}(R)$ the {\it embedding dimension} of $R$, i.e., the minimal number of generators of the maximal ideal $\m$. 

Let us begin with the following. 

\begin{lem}
For every parameter ideal $Q$ in $R$, one has an isomorphism
$$
\Ext_R^1(R/Q,R)\cong ((0):I)/Q
$$
of $R$-modules, where $I=(0): Q$. 
\end{lem}

\begin{proof}
Let $Q=(a)$. We first consider the following free resolution of $R/Q$
$$ 
\xymatrix{
\cdots  \ar[r] & R^l \ar[r]^{ \ [x_1 \cdots x_l]} & R \ar[r]^a & R \ar[r] & R/Q 
\ar[r] & 0, 
}
$$
where $l=\mu_R(I)$ and $I=(x_1, \dots , x_l)$. 
Taking the $R$-dual of this resolution, we have a complex
$$
\xymatrix{
0 \ar[r] & R \ar[r]^a & R \ar[r]^{^t[x_1 \cdots x_l]} 
& R^l \ar[r] & \cdots . 
}
$$
By this complex, we have an isomorphism $\Ext^1_R(R/Q, R) \cong ((0): I)/Q$. 
\end{proof}

Consequently, Question \ref{ques} is the same as the following.
$$ \mbox{Does it always hold that} \ (0): I \ne Q \ ? $$

We notice here that once $a$ is a non-zero divisor on $R$, then $\Ext_R^1(R/Q, R) \neq (0)$, because $(0):I=R \neq Q$. 
Hence $\Ext_R^1(R/Q, R) \neq (0)$ for every parameter ideal $Q= (a)$ in $R$, if $R$ is a Cohen-Macaulay local ring.

The following assertions are easy but we shall use them frequently in this paper. 

\begin{lem}\label{basic}
Let $Q$ be a parameter ideal in $R$ and $I=(0): Q$. Then we have 
the following. 
\begin{enumerate}
\item[{\rm (1)}] $Q \cap W = QW$.
\item[{\rm (2)}] $I \subseteq W$. 
\item[{\rm (3)}] $ \a \subseteq (0): I$. 
\item[{\rm (4)}] $Q: \m \subseteq QW: I$. 
\end{enumerate}
Furthermore, if the ideal $Q$ is standard, then we have
\begin{enumerate}
\item[$(5)$] $I=W$. 
\item[$(6)$] $\a = (0): I$. 
\item[$(7)$] $Q: \m \subseteq (0): I$. 
\end{enumerate}
\end{lem}

\begin{proof}
Let $Q=(a)$. Since $R/W$ is a $1$-dimensional Cohen-Macaulay local ring, 
the parameter $a$ is not a zero-divisor 
on $R/W$. Hence we have $Q \cap W = aW=QW$. $(2)$ follows from the fact 
$QI=(0)$. $(2)$ implies $(3)$. $(4)$ follows from the fact $I \subseteq \m$ and assertions $(1)$, $(2)$. Assume $QW=(0)$, that is, $W \subseteq I$. Then $(5)$ follows from $(2)$. $(5)$ implies $(6)$. Assertion $(7)$ follows from $(4)$. 
\end{proof}

\begin{prop}\label{std}
Let $Q$ be a parameter ideal in $R$.  
Then we have 
$$ \Ext_R^1(R/Q,R) \ne (0), $$ 
if either of the following conditions holds. 
\begin{enumerate}
\item[$(1)$] The ideal $Q$ is standard. 
\item[$(2)$] The ideal $Q^2$ is standard and $I=(0): Q$ is contained in $Q$. 
\end{enumerate}
\end{prop}

\begin{proof}
$(1)$ Suppose $\Ext_R^1(R/Q, R) = (0)$. 
Then we have the equality $Q:\m=Q$ by Lemma \ref{basic} $(7)$, which is impossible.

$(2)$ Let $Q=(a)$. We consider the following exact sequence.
$$ 
\xymatrix{
0 \ar[r] & R/(Q+I) \ar[r]^{\ \ \ a} & R/Q^2 \ar[r] & R/Q \ar[r] & 0. 
}
$$
Since $I \subseteq Q$, the above short exact sequence yields an exact sequence
$$
\xymatrix{
\Ext_R^1(R/Q, R) \ar[r] & \Ext_R^1(R/Q^2, R) \ar[r] & \Ext_R^1(R/Q, R). 
}
$$
Since the parameter ideal $Q^2$ is standard, 
we have $\Ext_R^1(R/Q^2, R) \ne (0)$ by $(1)$. 
Hence we get that $\Ext_R^1(R/Q, R) \ne (0)$. 
\end{proof}

It is well known that every parameter ideal $Q$ in a Buchsbaum local ring $R$ is standard. So, by Proposition 2.3 (1), we have the following. 

\begin{cor}
Assume that the ring $R$ is Buchsbaum. Then we have 
$$\Ext_R^1(R/Q,R)\ne (0)$$ 
for every parameter ideal $Q$ in $R$.
\end{cor}

\begin{thm}\label{thm a+z}
Suppose that a parameter ideal $Q=(a)$ is standard.
Then, for any $z \in W$, the element $a+z$ is a parameter of $R$, and 
$\Ext_R^1(R/(a+z),R) \ne (0)$. 
\end{thm}

\begin{proof}
Let $z \in W$. For every $\p \in \Min \, R$, since $z \in \p$, we have $a+z \notin \p$. Hence 
$a+z$ is a parameter for $R$. We put $b=a+z$. 
Then we have equalities
\begin{eqnarray*}
(0) & = & (b) [(0):b] \\
    & = & a [(0):b] + z[(0):b] \\
    & = & z [(0):b] \ \ (\mbox{since} \ (0):b \subseteq W \ \mbox{and} \ a \ \mbox{is standard}). 
\end{eqnarray*}
Therefore $z \in (0):((0):b)$. Suppose that $\Ext_R^1(R/(b),R) = (0)$. Then, since $\Ext_R^1(R/(b),R) \cong [(0):((0):b)]/(b)$ by Lemma 2.1, we have $z \in (0):((0):b)=(b)$. Hence we can write 
$z=by=(a+z)y$ for some $y \in R$. 
Then $y \in \m$ because $b \notin W$. Since $z(1-y)=ay$ and $1-y$ is a unit in $R$, we have $z \in Q \cap W =QW=(0)$. Hence $b=a$, which is a contradiction by Proposition \ref{std} (1). Therefore $\Ext_R^1(R/(b), R) \ne (0)$. 
\end{proof}

\begin{prop}\label{std2}
One has $\Ext_R^1(R/Q,R)\ne (0)$ for every parameter ideal $Q$ in $R$, 
if either of the following conditions holds. 
\begin{enumerate}
\item[$(1)$] The ideal $\a$ is not contained in $\m^2$. 
\item[$(2)$] $W^2=(0)$. 
\end{enumerate}
\end{prop}

\begin{proof}
Suppose $\Ext_R^1(R/Q, R) = (0)$ for some parameter ideal $Q$ in $R$. 

$(1)$ Take $x \in \a \backslash \m^2$. Since $\a \subseteq (0):I = Q$, 
we can write $x=ay$ for some $y \in R$. 
Then $y$ is a unit in $R$. Hence $\a = Q$. This implies that 
$Q$ is standard. By Proposition \ref{std} (1), this is impossible. 

$(2)$ By assumption $W^2=(0)$, $W \subseteq \a \subseteq (0):I=Q$. 
Hence $W \subseteq Q \cap W=QW$ and we have $W=(0)$ by Nakayama's lemma. 
Therefore $R$ is Cohen-Macaulay, which is a contradiction. 
\end{proof}

Before closing this section, let us give the following result. 

\begin{thm}
If $\mathrm{v}(R)\le 2$, then 
$\Ext_R^1(R/Q,R)\ne (0)$ for every parameter ideal $Q$ in $R$.
\end{thm}

\begin{proof}
We may assume that $\mathrm{v}(R)=2$. Furthermore, passing to the completion,  
we may assume that
$R$ is complete. Then there exists a 
two-dimensional regular local ring $S$ with the maximal ideal $\n$ 
such that $R \cong S/J$, where $J$ is 
an ideal in $S$ whose height is one. 
Since we may assume that $R$ is not Cohen-Macaulay, we can 
write $J=fL$ for some non-zero element $f \in \n$ and some $\n$-primary ideal 
$L$ in $S$. 
Since $W \cong (f)/J$ is the unmixed component of $R$, we have $W \cong S/L$. Therefore 
$LR \subseteq (0):W=\a$. 
Here, suppose that $\Ext_R^1(R/Q, R) =(0)$ for some parameter ideal $Q$ in $R$. 
Let $Q=gR$, where $g \in \n$. 
Then, since $LR \subseteq \a \subseteq (0):I=Q$, we have 
$L \subseteq (g)+J$ and hence $L = [L \cap (g)] + J$. Since 
$J \subseteq \n L$, $L = L \cap (g)$ by Nakayama's lemma. 
Hence we have $ L \subseteq (g)$. But this is 
impossible because $L$ is $\n$-primary and $\dim S =2$. 
\end{proof}



\section{Proof of Theorem \ref{main thm}}

The purpose of this section is to give a proof of Theorem \ref{main thm}. Recall that $R$ is a Noetherian local ring with the maximal ideal $\m$ 
and $\dim R=1$. Let $k=R/\m$ be the residue field of $R$. We set $W=H_{\m}^0(R)$ and $\a=(0):W$. We denote by $\mathrm{v}(R)$ the
embedding dimension of $R$. With these notation and assumption, 
we shall prove the following. 

\begin{thm}
Let $Q$ be a parameter ideal in $R$ and $I=(0):Q$. 
Suppose that $\m^2 W=(0)$. Then 
$$ \Ext_R^1(R/Q,R) \ne (0), $$
if one of the following holds.
\begin{enumerate}
\item[{\rm (1)}]
$\mathrm{v}(R)\le 4$.
\item[{\rm (2)}]
$\mu_R(W)\le 1$.
\item[{\rm (3)}]
$\mu_R(I)\le 1$.
\item[{\rm (4)}]
$\mathrm{v}(R/I)\le 2$.
\end{enumerate}
\end{thm}

\begin{proof}
Suppose $\Ext_R^1(R/Q, R) = (0)$. Then, since $\Ext_R^1(R/Q, R) \cong ((0):I)/Q$ by Lemma 2.1, we have $(0):I=Q$. 
Let $Q=(a)$. 
We first note that $\m^2 \subseteq \a$ and $Q^2$ is standard by the 
assumption $\m^2W=(0)$. 
By Proposition \ref{std2} (1), we have $\a \subseteq \m^2$ 
and hence $\a = \m^2$. Also, $a \notin \m^2$ because $Q$ is not standard by 
Proposition \ref{std} (1). Since $\m^2=\a \subset (0):I =Q$, we have 
$\m^2=a\m$. 
Note that $\m IW \subseteq \m^2 W=(0)$. Hence $\m W \subseteq ((0):I) \cap W=Q 
\cap W = QW$. Therefore $\m W=QW$. Furthermore one can check that $I$ is not contained in $\m^2$. 
Indeed, if $I \subseteq \m^2$, then $I \subseteq \m^2 = \a \subset Q$, which 
is impossible by Proposition \ref{std} (2). 

Now let $l=\ell_R((I+\m^2)/\m^2) >0$ and take $x_1, \dots , x_l \in I$ 
such that $\{ x_i \ \mod \ \m^2 \mid 1 \leq i \leq l \}$ is a $k$-basis of $(I+\m^2)/\m^2$. 
Then we have the following. 

\begin{claim}
\begin{enumerate}
\item[$({\rm{i}})$] $a, x_1, \dots , x_l$ is a part of 
a minimal system of generators for $\m$. 
\item[$({\rm{ii}})$] The equality $(0):(x_1, \dots , x_l)=(0):I$ holds. 
\end{enumerate}
\end{claim}

\begin{cpf}
$({\rm{i}})$ Let $\alpha, \beta_i \in R$ and suppose that 
$\alpha a + \sum_{i=1}^l \beta_i x_i \in \m^2$. 
Since $\m^2=a\m$, we can write $\alpha a + \sum_{i=1}^l \beta_i x_i = a \gamma$ for some $\gamma \in \m$. 
Hence $(\alpha - \gamma)a \in (x_1, \dots , x_l) \subseteq I=(0):a$. 
If $\alpha$ is a unit in $R$, then $\alpha - \gamma$ is also a unit in $R$ and hence $a \in I=(0):a$ and $a^2=0$, which is impossible because $a$ is a parameter. 
Therefore $\alpha \in \m$ and hence each $\beta_i \in \m$. 

$({\rm{ii}})$ $ I \subseteq (x_1, \dots , x_l)+\m^2$, so that the equality 
$I=[(x_1, \dots , x_l)+\m^2] \cap I=(x_1, \dots , x_l)+[\m^2 \cap I]$ holds. 
Since $\m^2 \cap I \subseteq Q \cap W = aW =\m W \subseteq (0):\m $, 
we then have the equality $(0):(x_1, \dots , x_l)=(0):I$. 
\qed
\end{cpf}

Let $\m=(a, x_1, \dots , x_l, x_{l+1}, \dots , x_n)$, where $n+1=\mathrm{v}(R)$. 
Since $\m^2=a\m$, we can write each $x_ix_j=a \delta_{ij}$ for some $\delta_{ij} \in \m$. 
Then we may assume that $\delta_{ij} \in 
(a, x_{l+1}, \dots , x_n)$ because $a x_i=0$ for all $1 \leq i \leq l$. 
Let $V$ be the $k$-subspace of $\m/\m^2$ spanned by 
$\{ \delta_{ij} \ \mod \ \m^2 \mid 1 \leq i \leq l,  \ 1 \leq j \leq n \}$ and let
$q=\dim_k V$. Then 

\begin{claim}
$q \leq n-l$. 
\end{claim}

\begin{cpf}
It is clear that $q+l \leq n+1$. Suppose $q+l=n+1$. 
Then $(a, x_{l+1}, \dots , x_n) \subseteq (\delta_{ij} \mid 1 \leq i \leq l, \ 1 \leq j \leq n)+\m^2$, so that $(a, x_{l+1}, \dots , x_n) = 
(\delta_{ij} \mid 1 \leq i \leq l, \ 1 \leq j \leq n)+[\m^2 \cap (a, x_{l+1}, \dots , x_n)]$. Since $\m^2=a\m$, we have $(a, x_{l+1}, \dots , x_n) = 
(\delta_{ij} \mid 1 \leq i \leq l, \ 1 \leq j \leq n)+\m(a, x_{l+1}, \dots , x_n)$. Hence $(a, x_{l+1}, \dots , x_n)=(\delta_{ij} \mid 1 \leq i \leq l, \ 1 \leq j \leq n)$ by Nakayama's lemma and hence the equality 
$$\m=(x_1, \dots , x_l)+(\delta_{ij} \mid 1 \leq i \leq l, \ 1 \leq j \leq n)
$$ holds. 
Let $\p \in \Min \, R $. Then $a \notin \p$. 
Since $ax_i=0$ for all $1 \leq i \leq l$, we have $(x_1, \dots , x_l) \subseteq \p$. Hence $a \delta_{ij} =x_i x_j \in \p$ for all $1 \leq i \leq l$ and $1 \leq j \leq n$. Therefore $(\delta_{ij} \mid 1 \leq i \leq l, \ 1 \leq j \leq n) 
\subseteq \p$. Hence $\m \subseteq \p$, which is impossible. Thus we have the inequality $q+l \leq n$. 
\qed
\end{cpf}

For any $n$-elements $a_1, \dots , a_n \in R$, we consider the following 
condition: 

\begin{eqnarray}\label{sharp}
c_i:=\sum_{j=1}^n a_j \delta_{ij} \in \m^2 \ \ \ \mbox{for all}  \ \ \ 1 \leq i \leq l. 
\end{eqnarray}

The elements $a_1, \dots , a_n \in R$ satisfying condition \ref{sharp} 
have the following property. 

\begin{claim}
If the elements $ a_1, \dots , a_n  \in R $ 
satisfy condition \ref{sharp}, then $a_i \in \m$ for all $1 \leq i \leq n$. 
\end{claim}

\begin{cpf}
For any $1 \leq i \leq l$, 
$$
ac_i  =  \sum_{j=1}^n a_j(a \delta_{ij}) 
      =  \sum_{j=1}^n a_j (x_i x_j) 
      =  x_i  \sum_{j=1}^n a_j x_j . 
$$
Hence we have 
$$a^2 c_i =ax_i \sum_{j=1}^n a_j x_j =0, $$
because $ax_i=0$ for all $1 \leq i \leq n$. 
Therefore each $c_i$ belongs to $W$. 
Since $a_1, \dots , a_n$ satisfy condition \ref{sharp}, 
$c_i \in \m^2 \cap W \subseteq Q \cap W = aW$. 
Hence $ac_i = x_i \left( \sum_{j=1}^n a_j x_j \right) =0 $ for all $1 \leq i \leq l$. 
Therefore we have 
$$\sum_{j=1}^n a_j x_j \in (0):(x_1, \dots , x_l)=(0):I=Q. $$
We write $\sum_{j=1}^n a_j x_j=az$ for some $z \in R$. 
Then $a_j \in \m$ for all $1 \leq j \leq n$, because the set $\{ a, x_1, \dots , x_n \}$ is a minimal system of generators for $\m$. 
\qed
\end{cpf}

\begin{claim}
$ql \geq n$. Hence we have $l \geq 2$ and $n-l \geq 2$. 
\end{claim}

\begin{cpf}
We take $\xi_1, \dots , \xi_q \in (\delta_{ij} \mid 1 \leq i \leq l, \ 
1 \leq j \leq n ) $ 
such that 
$ \{ \xi_k \ \mod \ \m^2  \mid  1 \leq k \leq q 
\} $ 
is a $k$-basis of $V$ and write 
$$\delta_{ij} \equiv \sum_{k=1}^q c^{ik}_{j} \xi_k \ \ \mod \ \m^2,  $$
where $c_j^{ik} \in R$. 
We now consider the following system of linear equations in variables $y_1, \dots , y_n$ over $k=R/\m$ :
\begin{equation}\label{star}
\sum_{j=1}^n \bar{c_j^{ik}} y_j =0 \ \ \ (1 \leq i \leq l, \  1 \leq k \leq q),  
\end{equation}
where $\bar{\ast}$ denotes the reduction of $\ast$ $\mod$ $\m$. 
Suppose $\bar{a_1}, \dots , \bar{a_n} \in k$ is a solution of \ref{star}. 
Then $\sum_{j=1}^n a_j c^{ik}_{j} \in \m$ for all 
$1 \leq i \leq l, \ 1 \leq k \leq q$. Therefore we have 
$$c_i = \sum_{j=1}^n a_j \delta_{ij} \equiv \sum_{j=1}^n a_j \left( \sum_{k=1}^q c^{ik}_{j} \xi_k \right) = \sum_{k=1}^q \left( \sum_{j=1}^n a_j c^{ik}_{j} 
\right) 
\xi_k \equiv 0 \ \mod \ \m^2. $$
Hence the elements $a_1, \dots , a_n$ satisfy condition \ref{sharp}. Thus $\bar{a_j}=0$ for all $1 \leq j \leq n$ by Claim 3. Therefore \ref{star} has the only trivial solution, which shows that $ql \geq n$. This implies $l \geq 2$ and $n-l \geq 2$. 
Indeed, if $l = 1$, then $n \leq q \leq n-1$ by Claim 2. This is impossible. 
Also, if $n-l \leq 1$, then $n \leq ql \leq (n-l)l \leq l \leq n$. Hence 
$n=ql=(n-l)l=l$. Therefore $q=1$ and $n=l$. Again, by Claim 2, this is impossible. 
\qed
\end{cpf}

%
%

Now suppose $\mathrm{v}(R) \leq 4$. Then $\mathrm{v}(R)-1=n=(n-l)+l \geq 4$, which is a contradiction. Suppose $\mu_R(W) \leq 1$. Since $(I+\m^2)/\m^2 \subseteq (W+\m^2)/\m^2=W/\m W$, it follows $l \leq 1$, which is a contradiction. 
If $\mu_R(I) \leq 1$, then $l \leq 1$. If $\mathrm{v}(R/I)=\ell_R(\m/(I+\m^2)) \leq 2$, 
then $l=\ell_R((I+\m^2)/\m^2) \geq \mathrm{v}(R)-2=(n+1)-2=n-1$. Hence $n-l \leq 1$, which is a contradiction. 

Consequently, we have $\Ext_R^1(R/Q, R) \ne (0)$. This is a proof of 
Theorem 3.1. 
\end{proof}

As a direct consequence, we have the following, which is Theorem \ref{main thm}.

\begin{cor}\label{main cor}
If $\m^2 W=(0)$ and $\mathrm{v}(R)\le 4$, we then have 
$$ \Ext_R^1(R/Q,R) \ne (0)  \ \mbox{for every parameter ideal} \ Q \ \mbox{in} \ R. $$
\end{cor}

\begin{rem}
The assumption in Corollary \ref{main cor} is the best possible. Indeed, 
in Section \ref{count ex}, we shall construct examples, showing that for a given 
integer $v \geq 5$, there exists a one-dimensional 
Noetherian local ring $(R, \m)$ 
with the embedding dimension $\mathrm{v}(R)=v$ and $\m^2 W=(0)$, which contains a parameter ideal $Q$ 
such that $\Ext^1_R(R/Q, R)=(0)$. 
\end{rem}

\section{Affirmative examples}

In this section, we shall give some affirmative examples. First, 
we give the following example, which follows from Theorem \ref{main thm}. 


\begin{ex}\label{but}
Let $k$ be a field and let $S=k[[X,Y,Z]]$ be a formal power series ring.
Set $U=(X,Y)$, $L=(X^2,XY-YZ,Y^2-XZ,Z^2)$ and $J=U \cap L$. We put $R=S/J$. Then $\dim R =1$ and we have $\Ext_R^1(R/Q, R) \ne (0)$ for every parameter ideal $Q$ in $R$. 
\end{ex}

\begin{proof}
One can check that $J=(X^2, XY-YZ, Y^2-XZ, XZ^2, YZ^2)$. Let $x, y$ denote respectively the reductions of $X, Y$ mod $J$. Then, since the unmixed component of $R$ is $\p:=(x, y)$, we have $W=\p$. It is easy to see that $\mathrm{v}(R)=3$ and $\m^2W=(0)$. Hence $\Ext_R^1(R/Q, R) \ne (0)$ for every parameter ideal $Q$ in $R$ by Theorem \ref{main thm}. 
\end{proof}

To construct another class of affirmative examples, we need the following. 

\begin{prop}\label{either}
Let $(S,\n)$ be a regular local ring with $\dim S>0$.
Let $U$, $L$ be ideals in $S$ satisfying the following three conditions. 
\begin{enumerate}
\item[{\rm (i)}]
The ring $S/U$ is a one-dimensional Cohen-Macaulay ring.
\item[{\rm (ii)}]
The ideal $L$ is an $\n$-primary ideal. 
\item[{\rm (iii)}]
The ideal $J:=U \cap L$ is contained in $\n^2$.
\end{enumerate}
We put $R=S/J$. Then $\dim R=1$ and we have 
$$\Ext_R^1(R/Q, R) \ne (0)  \ \mbox{for every parameter ideal} \ Q \ \mbox{in} \ 
R, $$
if either of the following conditions holds. 
\begin{enumerate}
\item[$(1)$] The ideal $L$ is not contained in $\n^2$. 
\item[$(2)$] The ideal $J=U \cap L$ is contained in $\n L$.
\end{enumerate}
\end{prop}

\begin{proof}
We may assume that the ideal $U$ is not contained in $L$. Since $W \cong U/J$, 
we have $LR \subseteq (0):W=\a$. 
Suppose $L$ is not contained in $\n^2$. 
Then $\a$ is not contained in $\m^2$, since $J=U \cap L \subseteq \n^2$. By Proposition 
\ref{std2} (1), we have that $\Ext_R^1(R/Q, R) \ne (0)$ 
for every parameter ideal $Q$ in $R$. 
Suppose $J \subseteq \n L$. Assume the contrary and choose a parameter ideal $Q = fR$ in $R$ such that $(0):I=Q$, where $f \in \n$. Since $\a \subseteq (0):I=Q$, we get $L \subseteq (f)+J$. Hence $L=[(f) \cap L] + \n L$. By Nakayama's lemma, we have $L=(f) \cap L$. Therefore $\dim S=1$ and hence $U=(0)$. This is a contradiction. 
\end{proof}

Using this, we have the following simple affirmative example, which 
does not follow from Theorem \ref{main thm}.  

\begin{ex}
Let $n>0$ and $m>l>0$ be integers. 
Let $k$ be a field and let $S=k[[X_1,X_2,\dots,X_n,Z]]$ be a formal power series ring. Set $U=(X_1^l, \dots, X_n^l)$, $L=(X_1^m, \dots, X_n^m, Z)$ and $J=U \cap L$. Then $J \subseteq \n^2$ and $L$ is not contained in $\n^2$, 
where $\n=(X_1,X_2,\dots,X_n,Z)$. Hence, for every parameter ideal $Q$ in 
$R:=S/J$, we have $\Ext_R^1(R/Q,R) \ne (0)$.
\end{ex}



The following example satisfies neither of the assumptions of Theorem \ref{main thm} and Proposition \ref{either}. But $\Ext_R^1(R/Q,R) \ne (0)$ holds for 
every parameter ideal $Q$ in $R$. 

\begin{ex}
Let $k$ be a field and let $S=k[[X,Y,Z]]$ be a formal power series ring.
Set $U=(X,Y)$, $L=(X^2,Y^2,Z^2)$ and $J=U \cap L$. We put $R=S/J$. Then $\dim R=1$ and we have $\Ext_R^1(R/Q, R) \ne (0)$ for every parameter ideal 
$Q$ in $R$. 
\end{ex}

\begin{proof} 
One can check that $J =(X^2,Y^2,XZ^2,YZ^2)$. Let $x, y, z$ denote respectively the reduction of $X, Y, Z$ mod $J$. Since the unmixed component of $R$ 
is $\p :=(x, y)$, we have $W=\p$. Then it is easy to see that $\m^2W \ne (0)$ and $\m^3W=(0)$. Suppose that there exists a parameter ideal $Q=(a)$ in $R$ such that 
$\Ext_R^1(R/Q, R) =(0)$. We may assume that $a=z^n+b$ where $n>0$ and 
$b \in \p$. Furthermore, since $z^2$ is standard, we may assume that $a=z+b$. 
Indeed, if $a=z^n+b$ for some $n \geq 2$, then $\Ext_R^1(R/(a), R) \neq (0)$ by Theorem \ref{thm a+z}, because $z^n$ is standard and $b \in W$. 

Since $xyW=(0)$, $xy \in (0):W = \a \subseteq (0):I=Q$. We write 
$xy=ac$ for some $c \in R$. Then $ c \in W$ because $ac \in W=\p$ and $a \notin \p$. Here we write 
\begin{eqnarray*}
b & = & b_1x+b_2y+b_3xy+b_4yz+b_5zx+b_6xyz, \\
c & = & c_1x+c_2y+c_3xy+c_4yz+c_5zx+c_6xyz,  
\end{eqnarray*}
where $b_i, c_i \in k$ for all $1 \leq i \leq 6$. 
Then we have 
\begin{eqnarray*}
xy &=& ac \\
   &=& (z+b)c \\
   &=& (b_1c_2+b_2c_1)xy+c_2yz+c_1zx+(c_3+b_1c_4+b_2c_5+b_4c_1+b_5c_2)xyz.  
\end{eqnarray*}
Since $\{xy, yz,zx,z^2 \}$ is a minimal system of generators 
for $\m^2$, we then have 
$c_1=c_2=0$ and $b_1c_2+b_2c_1=1$. This is a contradiction. 
\end{proof}

\section{Counterexamples}\label{count ex}

In this section, we will consider constructing examples which give a negative answer to Question \ref{ques}.

Let $n,l$ be integers with $2\le l\le n-2$.
Let $k$ be a field, $S=k[X_1,\dots,X_n,A]$ a polynomial ring.
The ring $S$ is a $\Z$-graded ring with $S_0=k$ and $\deg X_i = \deg A=1$ for $1\le i\le n$.
Set $V=\sum_{j=l+1}^n kX_j$.
Note that $\dim_kS_1=n+1\ge 5$.

\begin{lem}\label{delta}
There exists an $l\times n$ matrix $\Delta = (\Delta_{i,j})$ over $S$ which satisfies the following.
\begin{enumerate}
\item[{\rm (1)}]
The submatrix $(\Delta_{i,j})_{1\le i,j\le l}$ is symmetric.
\item[{\rm (2)}]
$V=\sum_{1\le i\le l,\,1\le j\le n}k\Delta_{i,j}$.
\item[{\rm (3)}]
If $c_1,\dots,c_n\in k$ satisfies $\Delta\left(
\begin{smallmatrix}
c_1\\
\vdots\\
\\
c_n
\end{smallmatrix}
\right)=0$, then $c_1=\cdots =c_n=0$.
\end{enumerate}
\end{lem}

\begin{pf}
If $l\le n-l$, then set 
$$
\Delta =
\begin{pmatrix}
X_{l+1} & X_{l+2} & \cdots & X_{2l} & 0 & \cdots & 0 \\
X_{l+2} & 0 & \cdots & 0 & 0 & \cdots & 0 \\
\vdots & \vdots & & \vdots & \vdots & & \vdots \\
X_{2l-1} & 0 & \cdots & 0 & 0 & \cdots & 0 \\
X_{2l} & 0 & \cdots & 0 & X_{l+1} & \cdots & X_n
\end{pmatrix}.
$$
It is easy to see that this matrix satisfies all the conditions in the lemma.

Let us consider the case where $l>n-l$.
Put $\alpha = n-l$.
We can write $l=\alpha q+r$ for some $0\le r<\alpha$.

\begin{claim*}
One has $0<q\le l-2$.
\end{claim*}

\begin{cpf}
If $q=0$, then $l=r<\alpha$, which is a contradiction.
Hence $q>0$.
Assume $q>l-2$.
Then $q\ge l-1$, and we have $l-r=\alpha q\ge \alpha(l-1)\ge 2(l-1)$ since $\alpha \ge 2$.
Hence $2\le l\le 2-r\le 2$.
Therefore we obtain $l=2$ and $r=0$.
It follows that $2=\alpha q$.
Since $\alpha\ge 2$, we get $\alpha =2=l>\alpha$.
This is a contradiction.
\qed
\end{cpf}

We construct a matrix $\Delta$ as follows:
$$
\begin{cases}
\Delta_{i,j}=\Delta_{j,i}=X_{l+j-\alpha(i-1)} & \text{ if }1\le i\le q\text{ and }\alpha(i-1)<j\le\alpha i, \\
\Delta_{q+1,j}=\Delta_{j,q+1}=X_{j+r} & \text{ if }\alpha q<j\le l, \\
\Delta_{l,j}=X_j & \text{ if }l<j, \\
\Delta_{i,j}=0 & \text{ otherwise}.
\end{cases}
$$
Then we can check that this matrix $\Delta$ satisfies the three conditions in the lemma.
\qed
\end{pf}

Let $\Delta$ be a matrix satisfying the conditions in Lemma \ref{delta}.
We define an ideal $J$ of $S$ as follows:
$$
J=(AX_1,\dots,AX_l)+(X_{l+1},\dots,X_n)^2+(X_iX_j-A\Delta_{i,j})_{1\le i\le l,\,1\le j\le n}.
$$
Since the submatrix $(\Delta_{i,j})_{1\le i,j\le l}$ is symmetric, we have
\begin{multline*}
J=(AX_1,\dots,AX_l)+(X_{l+1},\dots,X_n)^2\\
+(X_iX_j-A\Delta_{i,j})_{1\le i\le j\le l}+(X_iX_j-A\Delta_{i,j})_{1\le i\le l,\,l+1\le j\le n}.
\end{multline*}
Set $N=S_+=\bigoplus_{m>0}S_m\subseteq S$.
Lemma \ref{delta} says that $\Delta_{i,j}$ is in $V$.
Since $V$ is contained in $S_1$, the ideal $J$ is graded and contained in $N^2$.
Put $R=S/J$ and $M=R_+=N/J$.
Let $a,x_i,\delta_{i,j}$ be the residue classes of $A,x_i,\Delta_{i,j}$ in $R$, respectively.

\begin{prop}\label{prop}
One has the following.
\begin{enumerate}
\item[{\rm (1)}]
$\dim R=1$ and $\Min\,R=\{\p\}$, where $\p =(x_1,\dots,x_n)$.
\item[{\rm (2)}]
$M^2=aM$, $M^3=(a^3)$ and $M^2W=(0)$.
\item[{\rm (3)}]
$W=\p$ and $W_i=(0)$ for all $i\ge 3$.
\end{enumerate}
\end{prop}

\begin{pf}
(1) We make a claim.

\begin{claim*}
One has $\sqrt{J}=(X_1,\dots,X_n)$.
\end{claim*}

\begin{cpf}
Note that $(X_1,\dots,X_n)$ is a radical ideal of $S$.
Hence it is enough to show that $V(J)=V(X_1,\dots,X_n)$.
Let $P\in V(J)$.
Since $(X_{l+1},\dots,X_n)^2$ is contained in $J$, the elements $X_{l+1},\dots,X_n$ belong to $P$.
Hence $V$ is contained in $P$, and we have all $\Delta_{i,j}$ are in $P$ by Lemma \ref{delta}(2).
As $X_iX_j-A\Delta_{i,j}$ is in $J$, all $X_iX_j$ are in $P$.
In particular, $X_i^2$ is in $P$ for $1\le i\le l$.
Thus $X_i\in P$ for $1\le i\le l$.

Conversely, let $P\in V(X_1,\dots,X_n)$.
Then $V$ is contained in $P$, and Lemma \ref{delta}(2) says that all $\Delta_{i,j}$ are in $P$.
Hence all $X_iX_j-A\Delta_{i,j}$ are in $P$, and therefore $J$ is contained in $P$.
\qed
\end{cpf}

It follows from the above claim that $\dim R=\dim S/J=\dim S/\sqrt{J}=\dim k[A]=1$ and that $\min V(J)=\min V(\sqrt{J})=\min V(X_1,\dots,X_n)=\{(X_1,\dots,X_n)\}$, hence $\Min\,R=\{\p\}$.

(3) We begin with making the following claim.

\setcounter{claim}{0}
\begin{claim}\label{a}
One has $a^2x_i=0$ for $1\le i\le n$.
\end{claim}

\begin{cpf}
The claim is obvious for $1\le i\le l$, so let $l+1\le i\le n$.
Then $X_i$ is in $V=\sum_{1\le\alpha\le l,\,1\le\beta\le n}k\Delta_{\alpha,\beta}$, and we have $X_i=\sum_{\alpha,\beta}c_{\alpha,\beta}\Delta_{\alpha,\beta}$ for some $c_{\alpha,\beta}\in k$.
Hence $x_i=\sum_{\alpha,\beta}c_{\alpha,\beta}\delta_{\alpha,\beta}$, and we get $ax_i=\sum_{\alpha,\beta}c_{\alpha,\beta}(a\delta_{\alpha,\beta})=\sum_{\alpha,\beta}c_{\alpha,\beta}(x_\alpha x_\beta)$.
Therefore we obtain $a^2x_i = \sum_{\alpha,\beta}c_{\alpha,\beta}(ax_\alpha)x_\beta =0$ since $ax_\alpha =0$.
\qed
\end{cpf}

Note that $R/(a)$ is artinian.
Hence $a$ is a homogeneous parameter of $R$.
The above claim shows that $a^2\p=(0)$.
Since $(a)$ is $M$-primary, we have $M^s\p=(0)$ for some $s>0$.
It follows that $\p$ is contained in $W$.
On the other hand, as the ideal $W$ has finite length, it is nilpotent.
Therefore $W$ is contained in $\p$, and thus $W=\p$.

Here we make the following two claims:

\begin{claim}\label{i}
One has $ax_ix_j=0$ for $1\le i,j\le n$.
\end{claim}

\begin{cpf}
It holds that $ax_m=0$ if $1\le m\le l$.
Hence we may assume $l+1\le i,j\le n$.
Since $(x_{l+1},\dots,x_n)^2=(0)$, we have $x_ix_j=0$.
\qed
\end{cpf}

\begin{claim}\label{u}
One has $x_ix_jx_h=0$ for $1\le i,j,h\le n$.
\end{claim}

\begin{cpf}
The claim holds if $l+1\le i,j,h\le n$ since $(x_{l+1},\dots,x_n)^2=(0)$.
Let $1\le i\le l$.
Then $x_ix_j=a\delta_{i,j}$, and $x_ix_jx_h=a\delta_{i,j}x_h$.
Note that $\delta_{i,j}$ is in $kx_{l+1}+\cdots+kx_n$.
If $l+1\le h\le n$, then $\delta_{i,j}x_h=0$ as $(x_{l+1},\dots,x_n)^2=(0)$, and we get $x_ix_jx_h=0$.
If $1\le h\le l$, then $ax_h=0$, and $x_ix_jx_h=0$.
\qed
\end{cpf}

It follows from Claims \ref{a}, \ref{i} and \ref{u} that $W_i=(0)$ for all $i\ge 3$.

(2) Since $M=(x_1,\dots,x_n,a)=\p +(a)$, we have $M^2=aM+\p^2$.
For integers $\alpha,\beta$ with $1\le \alpha\le\beta\le n$, we have
$$
x_\alpha x_\beta =
\begin{cases}
a\delta_{\alpha,\beta} & (\alpha\le l), \\
0 & (l+1\le \alpha),
\end{cases}
$$
which is in $aM$.
Thus $M^2=aM$, and $M^3=aM^2=a^2M=a^2\p+(a^3)=(a^3)$ by Claim \ref{a}.
Added to it, we have $M^2W=M^2\p =aM\p=a(\p +(a))\p=a\p^2 = (0)$ by Claim 2.
\qed
\end{pf}

\begin{lem}\label{base}
The elements $ax_{l+1},\dots,ax_n,a^2$ form a $k$-basis of $R_2$.
\end{lem}

\begin{pf}
First of all, we claim the following.

\begin{claim*}
One has $J+(AX_{l+1},\dots,AX_n,A^2)=N^2$.
\end{claim*}

\begin{cpf}
Put $L=J+(AX_{l+1},\dots,AX_n,A^2)$.
It is obvious that $L$ is contained in $N^2$.
Fix integers $\alpha,\beta$ with $1\le \alpha \le \beta \le n$.
If $l+1\le \alpha$, then $X_\alpha X_\beta$ is in $J$, hence in $L$.
If $\alpha\le l$, then $X_\alpha X_\beta -A\Delta_{\alpha,\beta}$ is in $J$.
Since $\Delta_{\alpha,\beta}$ is in $V$, the element $A\Delta_{\alpha,\beta}$ is in the ideal $A(X_{l+1},\dots,X_n)$, which is contained in $L$.
Hence $X_\alpha X_\beta\in L$.
Thus the element $X_\alpha X_\beta$ is in $L$ for $1\le \alpha \le \beta \le n$.
Added to it, we have $NA=(AX_1,\dots, AX_l)+(AX_{l+1},\dots,AX_n,A^2)\subseteq L$.
Consequently, the ideal $N^2=(X_1,\dots,X_n)^2+NA$ is contained in $L$.
\qed
\end{cpf}

The above claim implies that $M^2=(ax_{l+1},\dots,ax_n,a^2)$.
Hence we have $R_2=(M^2)_2=k\cdot ax_{l+1}+\cdots +k\cdot ax_n +k\cdot a^2$.
Therefore $\dim_kR_2\le\alpha +1$, where $\alpha =n-l$.
Assume $\dim_kR_2<\alpha +1$.
Then $\dim_kS_2=\dim_kR_2+\dim_kJ_2<(\alpha+1)+\dim_kJ_2$.
We have 
\begin{multline*}
J=\underbrace{(AX_1,\dots,AX_l)}_{l} +\underbrace{(X_{l+1},\dots,X_n)^2}_{\frac{\alpha(\alpha+1)}{2}}\\
+\underbrace{(X_iX_j-A\Delta_{i,j})_{1\le i\le j\le l}}_{\frac{l(l+1)}{2}}+\underbrace{(X_iX_j-A\Delta_{i,j})_{1\le i\le l,\,l+1\le j\le n}}_{l\alpha}.
\end{multline*}
Hence $\dim_kJ_2\le l+\frac{\alpha(\alpha+1)}{2}+\frac{l(l+1)}{2}+l\alpha =\frac{n^2+n+2l}{2}$, and therefore $\dim_kS_2<(\alpha+1)+\frac{n^2+n+2l}{2}=\frac{(n+1)(n+2)}{2}=\dim_kS_2$.
This is a contradiction, and it must hold that $\dim_kR_2=\alpha+1$.
It follows that $ax_{l+1},\dots,ax_n,a^2$ form a $k$-basis of $R_2$.
\qed
\end{pf}

Now we are in the position to state and prove the main result of this section.

\begin{thm}
The element $a$ is a homogeneous parameter of $R$ satisfying $(0):((0):a)=(a)$, namely, $\Ext_R^1(R/(a),R)=(0)$.
\end{thm}

\begin{pf}
Set $I=(0):a$.
Let us prove the theorem step by step.

\setcounter{step}{0}
\begin{step}
The ideal $I$ is contained in $W$.

Indeed, since $(a)$ is an $M$-primary ideal, there is an integer $r>0$ such that $M^r$ is contained in $(a)$.
Since $aI=(0)$, we have $M^rI=(0)$.
\end{step}

\begin{step}
We have $I=I_1+I_2$.

Indeed, according to Proposition \ref{prop}(3), it holds that $W=W_0+W_1+W_2$.
Note that $W_0=W\cap R_0=W\cap k=(0)$.
Hence $W=W_1+W_2$.
Since $I$ is contained in $W$, we have $I=I_1+I_2$.
\end{step}

\begin{step}
We have $I_2=W_2$.

In fact, since $I$ is contained in $W$, $I_2$ is contained in $W_2$.
Proposition \ref{prop}(3) shows that $MW_2=(0)$, hence $aW_2=(0)$.
Thus $W_2$ is contained in $I=(0):a$.
Hence $W_2$ is contained in $I_2$.
\end{step}

\begin{step}
We have $I_2\subseteq (x_1,\dots,x_l)$.

In fact, since $I_2=W_2$, it suffices to check that $W_2$ is contained in $(x_1,\dots,x_l)$.
Let $\phi\in W_2$.
Note that $W_2=\p _2=\p\cap R_2=(x_1,\dots,x_n)\cap R_2$.
Hence we can write $\phi =\sum_{i=1}^nx_i\xi_i$ for some $\xi_i\in R_1$.
Let us show that $x_i\xi_i$ is in $(x_1,\dots,x_l)$ for $1\le i\le n$.
This is trivial if $1\le i\le l$, so let $l+1\le i\le n$.
Then $x_ix_j=0$ for $l+1\le j\le n$, and hence the element $x_ix_j$ is in $(x_1,\dots,x_l)$ for $1\le j\le n$.
As we saw in the proof of Claim \ref{a} in the proof of Proposition \ref{prop}, we can write $ax_i=\sum_{\alpha,\beta}c_{\alpha,\beta}(x_\alpha x_\beta)$ for some $c_{\alpha,\beta}\in k$.
Hence $ax_i\in (x_1,\dots,x_l)$.
Since $\xi_i$ belongs to $R_1=kx_1+\cdots+kx_n+ka$, we obtain $x_i\xi_i\in (x_1,\dots,x_l)$.
\end{step}

\begin{step}
We have $I=(x_1,\dots,x_l)$.

Indeed, as $ax_i=0$ for $1\le i\le l$, the ideal $(x_1,\dots,x_l)$ is contained in $I=(0):a$.
Suppose that $(x_1,\dots,x_l)$ is strictly contained in $I$, and choose a homogeneous element $\phi\in I-(x_1,\dots,x_l)$.
Since $I=I_1+I_2$, the element $\phi$ is in either $I_1$ or $I_2$.
However $I_2$ is contained in $(x_1,\dots,x_l)$, $\phi$ must be in $I_1$, hence in $W_1=\p_1=kx_1+\cdots+kx_n$.
Therefore $\phi=\psi+\sum_{i=l+1}^nc_ix_i$ for some $\psi\in kx_1+\cdots+kx_l$ and $c_i\in k$.
We have $0=a\phi=\sum_{i=l+1}^nc_i(ax_i)$ since $ax_j=0$ for $1\le j\le l$.
Lemma \ref{base} shows that $c_i=0$ for $l+1\le i\le n$, and $\phi =\psi\in (x_1,\dots,x_l)$.
This is a contradiction.
\end{step}

Now, we shall prove that $(0):((0):a)=(a)$.
It is trivial that $(0):((0):a)$ contains $(a)$.
Suppose that $(0):((0):a)=(0):I$ strictly contains the ideal $(a)$, and choose a homogeneous element $\phi\in ((0):I)-(a)$.
Then, since the ideal $M^2=aM$ is contained in $(a)$ and $\phi$ is not in $(a)$, $\phi$ is not in $M^2$.
Hence $\deg\phi\le 1$.
Assume that $\deg\phi =0$.
Then $\phi$ is in $k$ and is nonzero.
Since $\phi I=(0)$, we have $(0)=I=(x_1,\dots,x_l)$, which is a contradiction.
Thus $\deg\phi =1$, equivalently, the element $\phi$ is in $R_1$.
We can write $\phi=\sum_{j=1}^nc_jx_j+ca$ for some $c_j,c\in k$.
It holds that $(0)=\phi I=\phi (x_1,\dots,x_l)$, and $0=\phi x_i=\sum_{j=1}^n c_j(x_ix_j)+c(ax_i)=\sum_{j=1}^nc_j(a\delta_{i,j})=a\sum_{j=1}^nc_j\delta_{i,j}$ for $1\le i\le l$.
Hence $\sum_{j=1}^nc_j\delta_{i,j}\in ((0):a)=I=(x_1,\dots,x_l)$.
Noting that $\delta_{i,j}$ is in $V=kx_{l+1}+\cdots+kx_n$, we see that $\sum_jc_j\delta_{i,j}=0$ for $1\le i\le l$, and thus
$$
\Delta\left(
\begin{smallmatrix}
c_1\\
\vdots\\
\\
c_n
\end{smallmatrix}
\right)=0.
$$
By Lemma \ref{delta}(3) we have $c_i=0$ for $1\le i\le n$, and $\phi =ca\in (a)$, which is a contradiction.
This contradiction completes the proof of the theorem.
\qed
\end{pf}

\section{Modules of finite G-dimension}

In this section, we will consider a problem on modules of finite G-dimension.
We start by recalling the definition of G-dimension.

\begin{defn}
Let $R$ be a Noetherian ring.
\begin{enumerate}
\item[{\rm (1)}]
Let $(-)^\ast$ denote the $R$-dual functor $\Hom_R(-,R)$.
A finitely generated $R$-module $M$ is said to be {\it totally reflexive} if $M$ is isomorphic to $M^{\ast\ast}$ and $\Ext_R^i(M\oplus M^\ast,R)=(0)$ for all $i>0$.
\item[{\rm (2)}]
The {\it Gorenstein dimension} ({\it G-dimension} for short) of a nonzero $R$-module $M$, which is denoted by $\Gdim _R\,M$, is defined as the infimum of integers $r$ such that there exists an exact sequence
$$
0 \to X_r \to X_{r-1} \to \cdots \to X_0 \to M \to 0
$$
of $R$-modules, where each $X_i$ is totally reflexive.
The G-dimension of the zero module is defined as $-\infty$.
\end{enumerate}
\end{defn}

It is known that G-dimension has the following properties.
For the details, see \cite{ABr} and \cite{C}.

\begin{prop}\label{g}
Let $R$ be a Noetherian ring and $M$ a finitely generated $R$-module.
Then the following statements hold.
\begin{enumerate}
\item[{\rm (1)}]
There is an inequality $\Gdim_R\,M\le\pd_R\,M$.
\item[{\rm (2)}]
If $\Gdim_R\,M<\infty$, then $\Gdim_R\,M=\depth\,R-\depth _R\,M$.
\item[{\rm (3)}]
If $\Gdim_R\,M<\infty$, then $\Gdim_R\,M=\sup\{\, i\, |\, \Ext_R^i(M,R)\ne (0)\,\}$.
\end{enumerate}
\end{prop}

The third author gave the following conjecture in \cite{T}.

\begin{conj}
Let $R$ be a Noetherian local ring.
Suppose that there exists an $R$-module $M$ of finite length and finite G-dimension.
Then $R$ is Cohen-Macaulay.
\end{conj}

It is well-known that the statement with ``G-dimension'' replaced by ``projective dimension'' holds; it follows from the Peskine-Szpiro intersection theorem (cf. \cite[Proposition 6.2.4]{R}).

Let $R$ be a $d$-dimensional Noetherian local ring, and $M$ an $R$-module of finite length and finite G-dimension.
Then one has $\Gdim_R\,M=\depth\,R-\depth_R\,M=\depth\,R$ by Proposition \ref{g}(2), and hence $\Ext_R^i(M,R)=(0)$ for $i>\depth\,R$ by Proposition \ref{g}(3).
Therefore, if the $R$-module $M$ satisfies
$$
\Ext_R^d(M,R)\ne (0),
$$
then one must have $d\le\depth\,R$, that is to say, $R$ is Cohen-Macaulay.
So, if Question \ref{ques0} has an affirmative answer, then the above conjecture is true.
However, as we have already seen in the previous section, Question \ref{ques0} does not have an affirmative answer.

Now, we are interested in whether the example which we constructed in the previous section is a counterexample to the above conjecture or not.
The main result of this section is the following proposition, which says that it is not a counterexample.

\begin{prop}\label{pp}
Let $R$ be the ring and $a$ the homogeneous parameter of $R$ which are constructed in Section 5.
Then $R$ is a standard graded algebra over a field with $\dim R=1$ and $\depth\,R=0$ (hence $R$ is not Cohen-Macaulay) and $\Ext_R^1(R/(a),R)=(0)$, but the $R$-module $R/(a)$ is not of finite G-dimension.
\end{prop}

For a graded ring $R$ and a graded $R$-module $M$, we denote by $H_M(t)$ the Hilbert series of $M$.
To prove the above proposition, we prepare the following result, which is the main theorem in \cite{ABS}.

\begin{thm}[Avramov-Buchweitz-Sally]\label{abs}
Let $k$ be a field and $R$ a positively graded $k$-algebra.
Let $M,N$ be finitely generated graded $R$-modules with $\Ext_R^i(M,N)=(0)$ for $i\gg 0$.
Then
$$
\sum_i(-1)^i H_{\Ext_R^i(M,N)}(t)=\frac{H_M(t^{-1})\cdot H_N(t)}{H_R(t^{-1})}.
$$
\end{thm}

The lemma below follows from this theorem.

\begin{lem}
Let $R$ be a positively graded algebra over a field $k$.
Let $M$ be a graded totally reflexive $R$-module of finite length.
Then
$$
\ell_R(M)=\ell_R(M^\ast),
$$
where $(-)^\ast=\Hom_R(-,R)$.
\end{lem}

\begin{pf}
Since $\Ext_R^i(M,R)=(0)$ for $i>0$, Theorem \ref{abs} yields an equality
$$
H_{M^\ast}(t)=\frac{H_M(t^{-1})\cdot H_R(t)}{H_R(t^{-1})}.
$$
Note by definition that the dual module $M^\ast$ is also totally reflexive.
Replacing $M$ with $M^\ast$ in the above equality, we get
$$
H_M(t)=\frac{H_{M^\ast}(t^{-1})\cdot H_R(t)}{H_R(t^{-1})}.
$$
Thus we obtain the following two equalities:
$$
\begin{cases}
H_{M^\ast}(t)\cdot H_R(t^{-1})=H_M(t^{-1})\cdot H_R(t),\\
H_M(t)\cdot H_R(t^{-1})=H_{M^\ast}(t^{-1})\cdot H_R(t).
\end{cases}
$$
Therefore we obtain
\begin{equation}\label{hilb}
H_M(t)\cdot H_M(t^{-1})=H_{M^\ast}(t)\cdot H_{M^\ast}(t^{-1}).
\end{equation}
Since $M$ has finite length, we can write $H_M(t)=a_0+a_1t+\cdots+a_st^s$ for some integers $a_0,\dots,a_s$, and so $H_M(1)=a_0+a_1+\cdots +a_s=\ell_R(M)$.
Similarly we have $H_{M^\ast}(1)=\ell_R(M^\ast)$.
Substituting $t=1$ in the equality \eqref{hilb} yields $\ell_R(M)^2=\ell_R(M^\ast)^2$.
It follows that $\ell_R(M)=\ell_R(M^\ast)$, as desired.
\qed
\end{pf}

Now we can achieve the purpose of this section.

\begin{ppf}
Suppose that the $R$-module $R/(a)$ has finite G-dimension.
Then $\Gdim_R\,R/(a)=\depth\,R-\depth_R\,R/(a)=0$ since $\depth\,R=0$.
Hence $R/(a)$ is a totally reflexive $R$-module.

It is easy to see that $R/(a)$ is isomorphic to $k[X_1,\dots,X_n]/(X_1,\dots,X_n)^2$, which has dimension $n+1$ as a $k$-vector space.
Hence $\ell_R(R/(a))=n+1$.

On the other hand, the module $(R/(a))^\ast$ is isomorphic to the ideal $(0):a=I=(x_1,\dots,x_l)$, and it holds that $I=I_1+I_2$.
We have $I_1=k\cdot x_1+\cdots+k\cdot x_l$, hence $\dim_kI_1=l$.
The $k$-vector space $I_2$ is contained in $\sum_{1\le i\le l,\,1\le j\le n}k\cdot x_ix_j+\sum_{1\le i\le l}k\cdot ax_i$.
We have $x_ix_j=a\delta_{i,j}\in k\cdot ax_{l+1}+\cdots+k\cdot ax_n$ and $ax_i=0$, so $I_2$ is contained in $k\cdot ax_{l+1}+\cdots+ k\cdot ax_n$.
Conversely, since $a^2x_m=0$ for $l+1\le m\le n$, we get $ax_m\in I$.
It follows that $I_2=k\cdot ax_{l+1}+\cdots+k\cdot ax_n$.
Lemma \ref{base} guarantees that $ax_{l+1}, \dots,ax_n$ are linearly independent over $k$.
Therefore $\dim_kI_2=n-l$.
Consequently, we obtain equalities
$$
\ell_R((R/(a))^\ast)=\ell_R(I)=\dim_kI_1+\dim_kI_2=l+(n-l)=n.
$$
In particular, we get $\ell_R(R/(a))\ne \ell_R((R/(a))^\ast)$.
Theorem \ref{abs} gives a contradiction.
Thus, the $R$-module $R/(a)$ does not have finite G-dimension, and the proof is completed.
\qed
\end{ppf}



\end{document}